\documentclass[12pt]{article}
\usepackage{amstex,amssymb}
\begin{document}
%%%%%%%%%%%%%%%%%%%  Environments  %%%%%%%%%%%%%%%%%%%%%%%%%%
\newtheorem{thm}{Theorem}[section]
\newtheorem{lem}[thm]{Lemma}
\newtheorem{prop}[thm]{Proposition}
\newtheorem{conj}[thm]{Conjecture}
\newtheorem{cor}[thm]{Corollary}
\newenvironment{dfn}{\medskip\refstepcounter{thm}
\noindent{\bf Definition \thesection.\arabic{thm}\ }}{\medskip}
\newenvironment{ex}{\medskip\refstepcounter{thm}
\noindent{\bf Example \thesection.\arabic{thm}\ }}{\medskip}
\newenvironment{proof}{\medskip
\noindent{\it Proof.\ }}{\hfill\hskip 2em$\Box$\par}
%%%%%%%%%%%%%%%%%%%%%%  Macros  %%%%%%%%%%%%%%%%%%%%%%%%%%%%%
\def\bb{\mathbb}
\def\C{\mathbin{\bb C}}
\def\R{\mathbin{\bb R}}
\def\eq#1{{\rm(\ref{#1})}}
\def\d{{\rm d}}
\def\dim{\mathop{\rm dim}}
\def\Hol{\mathop{\rm Hol}}
\def\vol{\mathop{\rm vol}}
\def\ms#1{\vert#1\vert^2}
\def\nm#1{\Vert #1 \Vert}
\def\md#1{\vert #1 \vert}
\def\bnm#1{\bigl\Vert #1 \bigr\Vert}
\def\bmd#1{\bigl\vert #1 \bigr\vert}
%%%%%%%%%%%%%%%%%%%%%%%%%%%%%%%%%%%%%%%%%%%%%%%%%%%%%%%%%%%%%
%%%%%%%%%%%%%%%%%%%  Text of Paper  %%%%%%%%%%%%%%%%%%%%%%%%%
%%%%%%%%%%%%%%%%%%%%%%%%%%%%%%%%%%%%%%%%%%%%%%%%%%%%%%%%%%%%%
\title{Asymptotically Locally Euclidean \\ 
metrics with holonomy ${\rm SU}(m)$}
\author{Dominic Joyce \\
Lincoln College, Oxford, OX1 3DR, England}
\date{May 1999}
\maketitle

\section{Introduction}

Let $G$ be a finite subgroup of U$(m)$ acting freely on 
$\C^m\setminus\{0\}$. Then $\C^m/G$ has an {\it isolated 
quotient singularity} at 0. Suppose $(X,\pi)$ is a 
{\it resolution} of $\C^m/G$. Then $X$ is a noncompact 
complex manifold modelled at infinity on $\C^m/G$.

In this paper we will study K\"ahler metrics $g$ on $X$ which 
are also {\it Asymptotically Locally Euclidean}, or {\it ALE} 
for short. This means that $g$ approximates the Euclidean 
metric $h$ on $\C^m/G$ by $g=h+O(r^{-2m})$, with appropriate 
decay in the derivatives of $g$. We are particularly 
interested in {\it Ricci-flat} ALE K\"ahler manifolds. 

The {\it Calabi conjecture} \cite{Cala1} describes the possible 
Ricci curvatures of K\"ahler metrics on a fixed compact complex 
manifold $M$, in terms of the first Chern class $c_1(M)$ of $M$. 
It was proved by Yau \cite{Yau} in 1976. The following theorem 
is a corollary of Yau's proof.

\begin{thm} Let\/ $M$ be a compact complex manifold admitting 
K\"ahler metrics, with\/ $c_1(M)=0$. Then there is a unique 
Ricci-flat K\"ahler metric in each K\"ahler class on~$M$. 
\label{cyrfthm}
\end{thm}

Our main results are Theorems \ref{alerfthm} and \ref{alesumthm}. 
Theorem \ref{alerfthm} is an analogue of Theorem \ref{cyrfthm} 
for ALE K\"ahler manifolds. It says that if $X$ is a resolution 
of $\C^m/G$ with $c_1(X)=0$, that is, a {\it crepant resolution},
then every K\"ahler class of ALE K\"ahler metrics on $X$ contains 
a unique Ricci-flat K\"ahler metric. Theorem \ref{alesumthm} says 
that these metrics have holonomy ${\rm SU}(m)$. When $m=2$ the metrics 
were constructed explicitly by Kronheimer and others.

Section 2 defines ALE metrics and ALE K\"ahler metrics, 
and \S 3 states the main results of the paper, postponing
the proofs until \S 6, and gives some examples. Section
4 develops some analytical tools for ALE manifolds: Banach
spaces of functions called weighted H\"older spaces, and 
elliptic regularity theory for the Laplacian $\Delta$ on 
them. In \S 5 we discuss $k$-forms and de Rham cohomology 
on ALE manifolds.

Section 6 states a version of the Calabi conjecture 
for ALE manifolds. Only a sketch of the proof is given; a
complete proof, following Yau \cite{Yau}, will be given
in the author's book \cite[\S 8]{Joyc3}. We apply this 
Calabi conjecture to prove Theorem \ref{alerfthm}, and 
then prove Theorem~\ref{alesumthm}.

A number of other people have already written papers on 
noncompact versions of the Calabi conjecture, and I should 
at once admit that there is some overlap between their results 
and mine. In particular, Tian and Yau \cite{TiYa1,TiYa2} 
and independently Bando and Kobayashi \cite{BaKo1,BaKo2} prove 
the following result \cite[Cor.~1.1]{TiYa2}, \cite[Th.~1]{BaKo2}:
\medskip

\begin{thm} Let\/ $X$ be a compact K\"ahler manifold with\/ 
$c_1(X)>0$, and\/ $D$ a smooth reduced divisor on $X$ such that\/ 
$c_1(X)=\alpha[D]$ for some $\alpha>1$. Suppose $D$ admits a 
K\"ahler-Einstein metric with positive scalar curvature. Then 
$X\setminus D$ has a complete Ricci-flat K\"ahler metric.
\end{thm}

Also Tian and Yau give estimates on the decay of the curvature 
of their Ricci-flat metric. With a certain amount of work, the 
existence of the metrics of Theorem \ref{alerfthm} follows from 
the theorem above. But our estimates on the asymptotic behaviour 
of the metrics are stronger than those proved by Tian and Yau. 
For example, we show that the curvature is $O(r^{-2m-2})$ for 
large $r$, but Tian and Yau only show that it is $O(r^{-3})$, 
which is not good enough for the applications we have in mind.

In a sequel \cite{Joyc4} we will extend the material of this 
paper to construct a class of Ricci-flat K\"ahler metrics on 
crepant resolutions of {\it non-isolated} singularities 
$\C^m/G$, which we will call {\it Quasi-ALE metrics}. These 
metrics are not covered by the work of Tian and Yau or Bando 
and Kobayashi. 

The original motivation for this paper and \cite{Joyc4} 
is that ALE and Quasi-ALE metrics with holonomy SU(2), 
SU(3), SU(4) and Sp(2) are essential ingredients in 
a new construction by the author of compact manifolds with 
the exceptional holonomy groups $G_2$ and Spin(7), 
which generalizes that of \cite{Joyc1,Joyc2}. This 
construction will be described at length in the author's 
book \cite{Joyc3}, which also discusses the results of
this paper and~\cite{Joyc4}.

\section{Asymptotically Locally Euclidean metrics}

Suppose $G$ is a finite subgroup of ${\rm SO}(n)$ that acts freely 
on $\R^n\setminus\{0\}$. Then $\R^n/G$ has an {\it isolated 
quotient singularity} at 0. Let $h$ be the Euclidean metric on 
$\R^n$. Then $h$ is preserved by $G$, as $G\subset{\rm SO}(n)$, 
and so $h$ descends to $\R^n/G$. Let $r$ be the radius function 
on $\R^n/G$, that is, $r(x)$ is the distance from $0$ to $x$ 
calculated using $h$. We will define a natural class of noncompact 
Riemannian manifolds $(X,g)$ called {\it ALE manifolds}, that have 
one infinite end upon which the metric $g$ asymptotically resembles 
the metric $h$ on $\R^n/G$ for large~$r$.

\begin{dfn} Let $X$ be a noncompact manifold of dimension $n$, 
and $g$ a Riemannian metric on $X$. We say that $(X,g)$ is an 
{\it Asymptotically Locally Euclidean manifold} asymptotic to
$\R^n/G$, or an {\it ALE manifold} for short, and we say that 
$g$ is an {\it ALE metric}, if the following conditions hold.

There should exist a compact subset $S\subset X$ and a map
$\pi:X\setminus S\rightarrow\R^n/G$ that is a diffeomorphism 
between $X\setminus S$ and the subset $\{z\in\R^n/G:r(z)>R\}$ 
for some fixed $R>0$. Under this diffeomorphism, the push-forward
metric $\pi_*(g)$ should satisfy
\begin{equation}
\nabla^k\bigl(\pi_*(g)-h\bigr)=O(r^{-n-k})\quad
\text{on $\{z\in\R^n/G:r(z)>R\}$,} 
\label{aleeq}
\end{equation} 
for all $k\ge 0$. Here $\nabla$ is the Levi-Civita connection 
of $h$, and $T=O(r^{-j})$ if $\md{T}\le Kr^{-j}$ for some~$K>0$.
\label{aledef}
\end{dfn}

If $G=\{1\}$, so that $(X,g)$ is asymptotic to $\R^n$, then 
we call $(X,g)$ an {\it Asymptotically Euclidean manifold}, or 
{\it AE manifold}. We shall call the map $\pi:X\setminus S
\rightarrow\R^n/G$ an {\it asymptotic coordinate system} for $X$. 
Equation \eq{aleeq} says that towards infinity the metric $g$ on 
$X$ (and its derivatives) must converge to the Euclidean metric 
on $\R^n/G$, with a given rate of decay. We will explain 
in \S 3 why we have chosen the powers $r^{-n-k}$ here.

\begin{dfn} Let $(X,g)$ be an ALE manifold asymptotic 
to $\R^n/G$. We say that a smooth function 
$\rho:X\rightarrow[1,\infty)$ is a {\it radius function} 
on $X$ if, given any asymptotic coordinate system 
$\pi:X\setminus S\rightarrow\R^n/G$, we have
\begin{equation}
\nabla^k\bigl(\pi_*(\rho)-r\bigr)=O(r^{1-n-k})\quad
\text{on $\{z\in\R^n/G:r(z)>R\}$,} 
\end{equation}
for all $k\ge 0$. This condition is independent of the choice 
of asymptotic coordinate system, and radius functions exist 
for every ALE manifold.
\label{aleraddef}
\end{dfn}

A radius function is a function $\rho$ on $X$ that approximates 
the function $r$ on $\R^n/G$ near infinity. In doing analysis on 
ALE manifolds, we will find it useful to consider H\"older spaces 
of functions in which the norms are weighted by powers $\rho^\beta$ 
of a radius function. Note that by definition $\rho\ge 1$, so we 
do not have to worry about small values of~$\rho$.

Here is one way to think about ALE metrics. The manifold $X$ is 
noncompact, but it can be compactified in a natural way by adding the 
boundary ${\mathcal S}^{n-1}/G$ at infinity. So we can instead regard 
$X$ as a {\it compact manifold with boundary}. Then ALE metrics are 
metrics on $X$ satisfying a certain natural boundary condition. 

It is a general principle in differential geometry that most results 
about compact manifolds can also be extended to results about compact 
manifolds with boundary, provided the right boundary conditions are 
imposed in the problem. ALE manifolds are an example of this principle, 
because many results about compact Riemannian manifolds have natural 
analogues for ALE manifolds. 

Next we define {\it ALE K\"ahler metrics}. Suppose $G$ is a 
finite subgroup of U$(m)$ acting freely on $\C^m\setminus\{0\}$. 
Then $\C^m/G$ has an isolated quotient singularity at 0, 
and the standard Hermitian metric $h$ on $\C^m$ descends 
to $\C^m/G$. Let $r$ be the radius function on $\C^m/G$. 
Suppose $(X,\pi)$ is a {\it resolution} of $\C^m/G$, that is, 
$X$ is a normal nonsingular variety with a proper birational 
morphism $\pi:X\rightarrow\C^m/G$. Then we can consider metrics 
on $X$ which are both K\"ahler, and ALE.

\begin{dfn} Let $(X,\pi)$ be a resolution of $\C^m/G$, 
with complex structure $J$, and let $g$ be a K\"ahler metric on $X$. 
We say that $(X,J,g)$ is an {\it ALE K\"ahler manifold asymptotic to} 
$\C^m/G$, and that $g$ is an {\it ALE K\"ahler metric}, 
if for some $R>0$ we have
\begin{equation}
\nabla^k\bigl(\pi_*(g)-h\bigr)=O(r^{-2m-k})\quad
\text{on $\{z\in\C^m/G:r(z)>R\}$,} 
\label{alekeq}
\end{equation}
for all $k\ge 0$. We say that a smooth function 
$\rho:X\rightarrow[1,\infty)$ 
is a {\it radius function} on $X$ if $\rho=\pi^*(r)$ on the subset 
$\bigl\{x\in X:\pi^*(r)\ge 2\bigr\}$. A radius function exists for 
every ALE K\"ahler manifold.
\label{alekdef}
\end{dfn}

Because $X$ is a resolution of $\C^m/G$, it comes equipped with 
a resolving map $\pi:X\rightarrow\C^m/G$, which gives a natural 
asymptotic coordinate system for $X$. The consequence of using 
this preferred asymptotic coordinate system is that on an ALE 
K\"ahler manifold $(X,J,g)$, both the metric $g$ and the complex 
structure $J$ are simultaneously asymptotic to the metric and 
complex structure on $\C^m/G$. We also use $\pi$ to simplify 
the definition of radius function. 

In dimension 2 one can also desingularize $\C^2/G$ by {\it 
deformation}. By adopting a slightly more general definition of 
ALE K\"ahler manifold we can include deformations and resolutions 
of deformations of $\C^2/G$, and most of our results also apply to
them. This will be discussed in \cite[\S 8.9]{Joyc3}. However, 
by {\it Schlessinger's Rigidity Theorem} \cite{Schl}, if $m\ge 3$ 
then an isolated quotient singularity $\C^m/G$ admits no
nontrivial deformations.

\section{Ricci-flat ALE K\"ahler manifolds}

We now state some results on Ricci-flat ALE K\"ahler manifolds, 
and give some examples. The proofs will be deferred until \S 6. 
A resolution $(X,\pi)$ of $\C^m/G$ with $c_1(X)=0$ is called a 
{\it crepant resolution}, as in Reid \cite{Reid}. A great deal 
is known about the algebraic geometry of crepant resolutions, 
especially when $\dim X$ is 2 or 3. In particular, for $\C^m/G$ 
to admit a crepant resolution $G$ must be a subgroup of ${\rm SU}(m)$, 
and when $m$ is 2 or 3 a crepant resolution of $\C^m/G$ exists 
for every finite subgroup $G$ of~${\rm SU}(m)$.

Our first proposition shows that Ricci-flat ALE K\"ahler metrics 
exist only on {\it crepant resolutions}. The proof is elementary, 
and we omit it.

\begin{prop} Let\/ $G$ be a finite subgroup of\/ ${\rm U}(m)$ 
acting freely on $\C^m\setminus\{0\}$, let\/ $(X,\pi)$ be a 
resolution of\/ $\C^m/G$, and suppose $g$ is a Ricci-flat ALE 
K\"ahler metric on $X$. Then $X$ is a crepant resolution of\/ 
$\C^m/G$ and\/~$G\subset{\rm SU}(m)$.
\label{alerfcrepprop}
\end{prop}

Next we define {\it K\"ahler classes} and the {\it K\"ahler cone} 
for ALE manifolds.

\begin{dfn} Let $(X,J,g)$ be an ALE K\"ahler manifold asymptotic 
to $\C^m/G$ for some $m>1$, with K\"ahler form $\omega$. Then 
$\omega$ defines a de Rham cohomology class $[\omega]\in H^2(X,\R)$ 
called the {\it K\"ahler class} of $g$. Define the {\it K\"ahler 
cone} $\mathcal K$ of $X$ to be the set of K\"ahler classes 
$[\omega]\in H^2(X,\R)$ of ALE K\"ahler metrics on $(X,J)$. It 
is not difficult to prove that $\mathcal K$ is an open convex 
cone in $H^2(X,\R)$, which does not contain zero.
\end{dfn}

The following two theorems will be proved in~\S 6.

\begin{thm} Let\/ $G$ be a nontrivial finite subgroup of\/ 
${\rm SU}(m)$ acting freely on $\C^m\setminus\{0\}$, and\/ 
$(X,\pi)$ a crepant resolution of\/ $\C^m/G$. Then each 
K\"ahler class of ALE K\"ahler metrics on $X$ contains a 
unique Ricci-flat ALE K\"ahler metric $g$. The K\"ahler 
form $\omega$ of\/ $g$ satisfies
\begin{equation}
\pi_*(\omega)=\omega_0+A\,\d\d^c(r^{2-2m})+\d\d^c\chi
\label{pisomeq}
\end{equation}
on the set $\bigl\{z\in\C^m/G:r(z)>R\bigr\}$, where $A<0$ 
and\/ $R>0$ are constants, $\omega_0$ is the K\"ahler form 
of the Euclidean metric on $\C^m/G$, $r$ the radius 
function on $\C^m/G$, and\/ $\chi$ a smooth function on 
$\bigl\{z\in\C^m/G:r(z)>R\bigr\}$ such that\/ 
$\nabla^k\chi=O(r^{\gamma-k})$ for each\/ $k\ge 0$ 
and\/~$\gamma\in (1-2m,2-2m)$.
\label{alerfthm}
\end{thm} 

Theorem \ref{alerfthm} is the main result of this paper, 
and is an analogue of Theorem \ref{cyrfthm} for ALE 
K\"ahler manifolds. We use the notation that 
$\d^cf=i(\overline\partial-\partial)f$, when $f$ is a 
differentiable function on a complex manifold. Then $\d^c$ 
is a real operator, and~$\d\d^c=2i\partial\overline\partial$.

Note that because $A<0$ in Theorem \ref{alerfthm}, the term
$A\,\d\d^c(r^{2-2m})$ in \eq{pisomeq} is nonzero. Therefore
$\pi_*(g)-h$ decays with order exactly $O(r^{-2m})$, and 
similarly $\nabla^k(\pi_*(g)-h)$ decays with order exactly 
$O(r^{-2m-k})$. Thus in Definition \ref{aledef} the decay 
rates given in \eq{aleeq} are {\it sharp} for all Ricci-flat 
ALE K\"ahler metrics, and cannot be improved upon. This is 
why we chose the powers $r^{-n-k}$ in our definition \eq{aleeq} 
of ALE metrics.

\begin{thm} Let\/ $G$ be a nontrivial finite subgroup of\/ 
${\rm SU}(m)$ acting freely on $\C^m\setminus\{0\}$, let\/ 
$(X,\pi)$ be a crepant resolution of\/ $\C^m/G$, and let\/ 
$g$ be a Ricci-flat ALE K\"ahler metric on $X$. Then $g$ has 
holonomy~${\rm SU}(m)$.
\label{alesumthm}
\end{thm} 

For an introduction to holonomy groups of Riemannian manifolds, 
and the connection between Ricci-flat K\"ahler metrics and 
holonomy ${\rm SU}(m)$, see Salamon~\cite{Sala}.

\subsection{Examples}

ALE K\"ahler manifolds with holonomy SU(2) are already very well 
understood. Eguchi and Hanson \cite{EgHa} gave an explicit formula 
in coordinates for the metrics of ALE spaces with holonomy SU(2) 
asymptotic to $\C^2/\{\pm1\}$, and this was generalized by Gibbons 
and Hawking \cite{GiHa} to explicit expressions for ALE spaces 
asymptotic to $\C^2/{\bb Z}_k$ for $k\ge 2$. More generally, 
Kronheimer \cite{Kron1,Kron2} gave an explicit, algebraic 
construction of every ALE manifold with holonomy SU(2), using 
the hyperk\"ahler quotient.

Thus, we can write down many explicit examples of ALE manifolds 
with holonomy SU(2). For $m\ge 3$, Calabi \cite[p.~285]{Cala2}
found an explicit ALE K\"ahler manifold with holonomy ${\rm SU}(m)$ 
asymptotic to $\C^m/{\bb Z}_m$, which we describe next. In the case 
$m=2$, Calabi's example coincides with the Eguchi--Hanson metric.

\begin{ex} Let $\C^m$ have complex coordinates $(z_1,\ldots,z_m)$,
let $\zeta={\rm e}^{2\pi i/m}$, and let $\alpha$ act on $\C^m$ by
$\alpha:(z_1,\ldots,z_m)\mapsto(\zeta z_1,\ldots,\zeta z_m)$. Then
$\alpha^m=1$, and the group $G=\langle\alpha\rangle$ generated by 
$\alpha$ is a subgroup of ${\rm SU}(m)$ isomorphic to ${\bb Z}_m$, 
which acts freely on $\C^m\setminus\{0\}$. Thus the quotient 
$\C^m/G$ has an isolated singular point at 0. Let $(X,\pi)$ be the 
{\it blow-up} of $\C^m/G$ at 0, so that $\pi^{-1}(0)\cong
\bb{CP}^{m-1}$. It is easy to show that $X$ is in fact a 
{\it crepant resolution} of~$\C^m/G$. 

Let $r$ be the radius function on $\C^m/G$, and define
$f:\C^m/G\,\big\backslash\{0\}\rightarrow\R$ by
\begin{equation}
f=\sqrt[m]{r^{2m}+1}+{1\over m}\sum_{j=0}^{m-1}\zeta^j
\log\left(\sqrt[m]{r^{2m}+1}-\zeta^j\right).
\label{calafdef}
\end{equation}
To define the logarithm of the complex number 
$\sqrt[m]{r^{2m}+1}-\zeta^j$ we cut $\C$ along the negative 
real axis, and set $\log(R{\rm e}^{i\theta})=\log R+i\theta$ 
for $R>0$ and $\theta\in(-\pi,\pi)$. Then $f$ is well-defined, 
and it is a smooth {\it real}\/ function on 
$\C^m/G\,\big\backslash\{0\}$, despite its complex definition. 

Define a (1,1)-form $\omega$ on $X\setminus\pi^{-1}(0)$ by 
$\omega=\d\d^c\pi^*(f)$. It can be shown that $\omega$ extends 
to a smooth, closed, positive (1,1)-form on all of $X$. Let $g$ 
be the K\"ahler metric on $X$ with K\"ahler form $\omega$. Then 
Calabi \cite[\S 4]{Cala2} shows that $g$ is complete and Ricci-flat, 
with $\Hol(g)={\rm SU}(m)$. Equation \eq{calafdef} is derived from 
\cite[eqn~(4.14), p.~285]{Cala2}. Note also that the action of 
${\rm U}(m)$ on $\C^m$ pushes down to $\C^m/G$ and lifts through 
$\pi$ to $X$, and $g$ is invariant under this action of 
${\rm U}(m)$ on~$X$. 

This metric $g$ on $X$ is an {\it ALE K\"ahler metric}. To 
prove this, we show using \eq{calafdef} that
\begin{equation}
f=r^2-{1\over m(m\!-\!1)}\,r^{2-2m}+O(r^{-2m})\quad
\text{on $\C^m/G\,\big\backslash\{0\}$, for large $r$.}
\end{equation}
Now the K\"ahler form of the Euclidean metric on $\C^m/G$ is 
$\omega_0=\d\d^c(r^2)$. Hence
\begin{equation}
\pi_*(\omega)=\omega_0-\frac{1}{m(m\!-\!1)}\,\d\d^c(r^{2-2m})
+\d\d^c\chi \quad\text{on $\C^m/G\,\big\backslash\{0\}$,}
\label{pisomeqex}
\end{equation}
where $\chi=f-r^2+{1\over m(m-1)}r^{2-2m}$. It is easy to show 
that $\nabla^k\chi=O(r^{-k-2m})$ on $\C^m/G\,\big\backslash\{0\}$
for large $r$, and it quickly follows that $g$ is an ALE K\"ahler 
metric on $X$, by Definition \ref{alekdef}. Also, $g$ is one 
of the Ricci-flat ALE K\"ahler metrics of Theorem \ref{alerfthm}, 
and comparing \eq{pisomeqex} with \eq{pisomeq} we see that 
$A=-\,{1\over m(m-1)}$, which verifies that~$A<0$.
\label{alecalabiex}
\end{ex}

For $m\ge 3$, the metrics of Example \ref{alecalabiex} are the
{\it only} explicit examples of ALE metrics with holonomy 
${\rm SU}(m)$ that are known, at least to the author. It is 
possible to find these metrics explicitly because they have a 
large symmetry group ${\rm U}(m)$, whose orbits are of real 
codimension 1 in $X$. Because of this, the problem can be 
reduced to a nonlinear, second-order ODE in one real variable, 
which can then be explicitly solved.

It is a natural question whether we can find an explicit, 
algebraic form for any or all of the other ALE metrics with 
holonomy ${\rm SU}(m)$ for $m\ge 3$, that exist on crepant 
resolutions of $\C^m/G$ by Theorem \ref{alerfthm}. The author 
believes that general ALE metrics with holonomy ${\rm SU}(m)$ 
for $m\ge 3$ are essentially transcendental, nonalgebraic 
objects, and that one cannot write them down explicitly using 
simple functions. Furthermore, the author conjectures that for 
$m\ge 3$, the metrics of Example \ref{alecalabiex} are the only 
ALE metrics with holonomy ${\rm SU}(m)$ that can be written 
down explicitly in coordinates.

\section{Analysis on ALE manifolds}

Let $(M,g)$ be a Riemannian manifold. Then the {\it H\"older 
spaces} $C^{k,\alpha}(M)$ are Banach spaces of functions on 
$M$, defined in Besse \cite[p.~456-7]{Bess}. When $M$ is 
compact, elliptic operators such as the Laplacian $\Delta$ 
have very good regularity properties on H\"older spaces. 
Here is a typical elliptic regularity result, following 
from \cite[Th.~27 \& Th.~31, p.~463-4]{Bess}. Theorems of 
this kind are essential tools in analytic problems such 
as the proof of the Calabi conjecture.

\begin{thm} Let\/ $(M,g)$ be a compact Riemannian manifold, 
let\/ $k\ge 0$ be an integer, and\/ $\alpha\in (0,1)$. Then 
for each\/ $f\in C^{k,\alpha}(M)$ with\/ $\int_Mf\,\d V_g=0$
there exists a unique $u\in C^{k+2,\alpha}(M)$ with\/ 
$\int_Mu\,\d V_g=0$ and\/ $\Delta u=f$. Moreover,
$\nm{u}_{C^{k+2,\alpha}}\le C\nm{f}_{C^{k,\alpha}}$ for some 
$C>0$ independent of\/ $u$ and\/~$f$. 
\label{mschauderthm}
\end{thm}

However, if $(X,g)$ is an ALE manifold then the results of 
Theorem \ref{mschauderthm} are false for $X$. This tells us 
that the $C^{k,\alpha}(X)$ are not good choices of Banach 
spaces of functions for studying elliptic operators on an ALE 
manifold. Instead, it turns out to be helpful to introduce 
{\it weighted H\"older spaces}, which we define next.

\begin{dfn} Let $(X,g)$ be an ALE manifold asymptotic to $\R^n/G$, 
and $\rho$ a radius function on $X$. For $\beta\in\R$ and $k$ a 
nonnegative integer, define $C^k_\beta(X)$ to be the space of 
continuous functions $f$ on $X$ with $k$ continuous derivatives, 
such that $\rho^{j-\beta}\bmd{\nabla^jf}$ is bounded on $X$ for 
$j=0,\ldots,k$. Define the norm $\nm{\,.\,}_{\smash{C^k_\beta}}$ on 
$C^k_\beta(X)$ by 
\begin{equation}
\nm{f}_{C^k_\beta}=\sum_{j=0}^k\sup_X\bmd{\rho^{j-\beta}\nabla^jf}.
\end{equation}
Let $\delta(g)$ be the injectivity radius of $g$, and write $d(x,y)$ 
for the distance between $x,y$ in $X$. For $T$ a tensor field on $X$ 
and $\alpha,\gamma\in\R$, define
\begin{equation}
\bigl[T\bigr]_{\alpha,\gamma}=
\sup\begin{Sb}x\ne y\in X\\d(x,y)<\delta(g)\end{Sb}
\left[\min\bigl(\rho(x),\rho(y)\bigr)^{-\gamma}\cdot
{\bmd{T(x)-T(y)}\over d(x,y)^\alpha}\right].
\label{aleholdereq}
\end{equation}
Here we interpret $\md{T(x)-T(y)}$ using parallel translation
along the unique geodesic of length $d(x,y)$ joining $x$ and~$y$.

For $\beta\in\R$, $k$ a nonnegative integer, and $\alpha\in(0,1)$, 
define the {\it weighted H\"older space} $C^{k,\alpha}_\beta(X)$ 
to be the set of $f\in C^k_\beta(X)$ for which the norm
\begin{equation}
\bnm{f}_{C^{k,\alpha}_\beta}=\bnm{f}_{C^k_\beta}
+\bigl[\nabla^k f\bigr]_{\alpha,\beta-k-\alpha}
\label{aleholdnormeq}
\end{equation}
is finite. Define $C^\infty_\beta(X)$ to be the intersection 
of the $C^k_\beta(X)$ for all $k\ge 0$. Both $C^k_\beta(X)$ 
and $C^{k,\alpha}_\beta(X)$ are Banach spaces, but 
$C^\infty_\beta(X)$ is not a Banach space.
\label{wholderdef}
\end{dfn}

This definition is taken from Lee and Parker \cite[\S 9]{LePa}.
A function $f$ in $C^k_\beta(X)$ or $C^{k,\alpha}_\beta(X)$ 
grows at most like $\rho^\beta$ as $\rho\rightarrow\infty$, 
and so the index $\beta$ should be interpreted as an 
{\it order of growth}. Similarly, the derivatives 
$\nabla^jf$ grow at most like $\rho^{\beta-j}$ for 
$j=1,\ldots,k$. As vector spaces of functions $C^k_\beta(X)$ 
and $C^{k,\alpha}_\beta(X)$ are independent of the choice of 
radius function $\rho$. The norms on these spaces do depend 
on $\rho$, but not in a significant way, as all choices of 
$\rho$ give equivalent norms. 

There is also another useful class of Banach spaces on ALE 
manifolds, the {\it weighted Sobolev spaces} $L^q_{k,\beta}(X)$,
which we will not define. They have similar analytic properties 
to the weighted H\"older spaces, and are described in 
\cite[\S 9]{LePa}. We have chosen to use weighted H\"older 
spaces instead, as they are often more convenient for 
nonlinear problems.

Next we discuss the analysis of the Laplacian $\Delta$ on ALE
manifolds. Much work has been done on the behaviour of $\Delta$ 
on weighted Sobolev spaces and H\"older spaces on $\R^n$, and 
more generally on AE manifolds. A useful guide, with references, 
can be found in Lee and Parker \cite[\S 9]{LePa}. Most of these
results apply immediately to ALE manifolds, with only very 
minor cosmetic changes to their proofs. 

\begin{prop} Let\/ $(X,g)$ be an ALE manifold of dimension 
$n$ asymptotic to $\R^n/G$, let\/ $\beta,\gamma\in\R$ satisfy 
$\beta+\gamma<2-n$, and suppose $u\in C^2_\beta(X)$ and\/ 
$v\in C^2_\gamma(X)$. Then
\begin{equation}
\int_X u\,\Delta v\,\d V_g=\int_X v\,\Delta u\,\d V_g.
\label{intxuveq}
\end{equation}
Let\/ $\rho$ be a radius function on $X$. Then 
$\Delta(\rho^{2-n})\in C^\infty_{-2n}(X)$ and
\begin{equation}
\int_X\Delta(\rho^{2-n})\,\d V_g=
{(n\!-\!2)\,\Omega_{n-1}\over\md{G}},
\label{intderhoeq}
\end{equation}
where $\Omega_{n-1}$ is the volume of the unit sphere 
${\mathcal S}^{n-1}$ in~$\R^n$.
\label{alede1prop}
\end{prop}

\begin{proof} Let $S_R$ be the subset $\{x\in X:\rho(x)\le R\}$ 
in $X$. Stokes' Theorem gives that
\begin{equation}
\int_{S_R}\bigl(u\Delta v-v\Delta u\bigr)\,\d V_g=\int_{\partial S_R}
\bigl[(u\nabla v-v\nabla u)\cdot{\bf n}\bigr]\,\d V_g, 
\label{uvsreq}
\end{equation}
where $\bf n$ is the inward-pointing unit normal to $\partial S_R$.
But for large $R$ we have $\vol(\partial S_R)=O(R^{n-1})$ and
$u\nabla v-v\nabla u=O(R^{\beta+\gamma-1})$ on $\partial S_R$, 
so that the r.h.s.~of \eq{uvsreq} is $O(R^{\beta+\gamma+n-2})$. 
Since $\beta+\gamma<2-n$ we see that the r.h.s.~of \eq{uvsreq} tends 
to zero as $R\rightarrow\infty$, and this proves~\eq{intxuveq}.

The point about the power $\rho^{2-n}$ is that $\Delta(r^{2-n})=0$ 
away from 0 in $\R^n/G$. Using the definitions of radius function 
and ALE manifold one can show that $\Delta(\rho^{2-n})\in 
C^\infty_{-2n}(X)$, as we want. Using Stokes' Theorem again 
we find that
\begin{equation}
\int_{S_R}\Delta(\rho^{2-n})\,\d V_g=\int_{\partial S_R}
\bigl[\nabla(\rho^{2-n})\cdot{\bf n}\bigr]\,\d V_g.
\end{equation}
But for large $R$ we have $\nabla(\rho^{2-n})\cdot{\bf n}\approx
(n\!-\!2)R^{1-n}$ and $\vol(S_R)\approx R^{n-1}\Omega_{n-1}/\md{G}$. 
Thus, letting $R\rightarrow\infty$ gives~\eq{intderhoeq}.
\end{proof}

\begin{thm} Let\/ $n>2$ and\/ $k\ge 0$ be integers and\/ 
$\alpha\in(0,1)$, and let\/ $\R^n$ have its Euclidean metric. Then
\begin{itemize}
\item[{\rm(a)}] Suppose $\beta\in(-n,-2)$. Then for each\/ 
$f\in C^{k,\alpha}_\beta(\R^n)$ there is a unique 
$u\in C^{k+2,\alpha}_{\beta+2}(\R^n)$ with\/~$\Delta u=f$.
\item[{\rm(b)}] Suppose 
$\beta\in(-1-n,-n)$. Then for each\/ 
$f\in C^{k,\alpha}_\beta(\R^n)$ there exists
$u\in C^{k+2,\alpha}_{\beta+2}(\R^n)$ with\/ $\Delta u=f$ if 
and only if\/ $\int_{\R^n}f\,\d V=0$, and $u$ is then unique.
\end{itemize}
In each case $\nm{u}_{\smash{C^{k+2,\alpha}_{\beta+2}}}\le 
C\nm{f}_{\smash{C^{k,\alpha}_\beta}}$ for some $C>0$ depending 
only on $n,k,\alpha$ and\/~$\beta$.
\label{alede1thm}
\end{thm}

\begin{proof} This is an analogue for $\R^n$ of Theorem 
\ref{mschauderthm}. If $u\in C^2_{\beta+2}(\R^n)$ for $\beta<-2$ 
and $\Delta u=f$, then by \cite[\S 2.4]{GiTr} we have
\begin{equation}
u(y)={1\over(n-2)\Omega_{n-1}}\int_{x\in\R^n}\md{x-y}^{2-n}f(x)\d x,
\label{intufeq}
\end{equation}
where $\Omega_{n-1}$ is the volume of the unit sphere 
${\mathcal S}^{n-1}$ in $\R^n$. This is {\it Green's 
representation} for $u$. Let $\rho$ be a radius function on $\R^n$. 
Then $\md{f(x)}\le\nm{f}_{\smash{C^0_\beta}}\rho(x)^\beta$, so 
\eq{intufeq} gives
\begin{equation}
\md{u(y)}\le{1\over(n-2)\Omega_{n-1}}\,\nm{f}_{C^0_\beta}
\cdot\int_{x\in\R^n}\md{x-y}^{2-n}\rho(x)^\beta\d x.
\end{equation}

We split this into integrals over the three regions 
$\md{x}\le{1\over 2}\md{y}$, ${1\over 2}\md{y}<\md{x}\le 2\md{y}$ 
and $\md{x}>2\md{y}$ in $\R^n$. Estimating the integral on each 
region separately we prove
\begin{equation}
\int_{x\in\R^n}\md{x-y}^{2-n}\rho(x)^\beta\d x\le
\begin{cases}
C'\rho(y)^{\beta+2} & \text{for $\beta\in(-n,-2)$,} \\
C'\rho(y)^{2-n} & \text{for $\beta<-n$.}
\end{cases}
\end{equation}
In case (a), if $\beta\in(-n,-2)$ then 
$\md{u(y)}\le C''\nm{f}_{\smash{C^0_\beta}}\rho(y)^{\beta+2}$ 
for some $C''>0$ depending only on $n$ and $\beta$, and so 
$u\in C^0_{\beta+2}(\R^n)$ and~$\nm{u}_{\smash{C^0_{\beta+2}}}
\le C''\nm{f}_{\smash{C^0_\beta}}$. 

One can extend this to show that 
$u\in C^{k+2,\alpha}_{\beta+2}(\R^n)$ and 
$\nm{u}_{\smash{C^{k+2,\alpha}_{\beta+2}}}\le 
C\nm{f}_{\smash{C^{k,\alpha}_\beta}}$ for some $C>0$ using the
method of {\it Schauder estimates}, as in \cite[\S 6]{GiTr}.
The difficulty in doing this is to correctly include the powers 
of $\rho$ involved in the weighted H\"older norm. To do this, 
for each $x\in\R^n$ we consider the ball $B_{\rho(x)/2}(x)$ 
of radius ${1\over 2}\rho(x)$ about $x$ in~$\R^n$.

On this ball we have $u=O\bigl(\rho(x)^{\beta+2}\bigr)$,
$\nabla^jf=O\bigl(\rho(x)^{\beta-j}\bigr)$ for $j=0,\dots,k$, and 
$[\nabla^kf]_\alpha=O\bigl(\rho(x)^{\beta-k-\alpha}\bigr)$. Using 
the Schauder interior estimates on the unit ball in $\R^n$ and 
rescaling distances by a factor ${1\over 2}\rho(x)$, we show that 
$\nabla^ju=O\bigl(\rho(x)^{\beta+2-j}\bigr)$ for $j=0,\dots,k\!+\!2$ 
and $[\nabla^{k+2}u]_\alpha=O\bigl(\rho(x)^{\beta-k-\alpha}\bigr)$
on the interior of $B_{\rho(x)/2}(x)$. Thus 
$u\in C^{k+2,\alpha}_{\beta+2}(\R^n)$ and 
$\nm{u}_{\smash{C^{k+2,\alpha}_{\beta+2}}}\le 
C\nm{f}_{\smash{C^{k,\alpha}_\beta}}$, completing the 
proof of case~(a).

Next we prove (b). Suppose $\beta\in(-1-n,-n)$,
$u\in C^{k+2,\alpha}_{\beta+2}(\R^n)$ and $\Delta u=f$. Then 
\begin{equation}
\int_{\R^n}f\,\d V=\int_{\R^n}1\Delta u\,\d V
=\int_{\R^n}u\Delta(1)\,\d V=0
\end{equation}
by Proposition \ref{alede1prop}, since $u\in C^2_{\beta+2}(\R^n)$ 
and $1\in C^2_0(\R^n)$ and $\beta+2+0<2-n$. Thus, given 
$f\in C^{k,\alpha}_\beta(\R^n)$, there can only exist 
$u\in C^{k+2,\alpha}_{\beta+2}(\R^n)$ with $\Delta u=f$ if 
$\int_{\R^n}f\,\d V=0$. So suppose that $\int_{\R^n}f\,\d V=0$, 
and define $u$ by
\begin{equation}
u(y)={1\over(n-2)\Omega_{n-1}}\int_{x\in\R^n}
\Bigl[\md{x-y}^{2-n}-\rho(y)^{2-n}\Bigr]f(x)\d x.
\label{intufreq}
\end{equation}
Since $\int_{\R^n}f\,\d V=0$ the term involving $\rho(y)^{2-n}$
in this integral vanishes, so the equation reduces to \eq{intufeq}
and thus $\Delta u=f$. From \eq{intufreq} we see that
\begin{equation*}
\md{u(y)}\le{1\over(n-2)\Omega_{n-1}}\,\nm{f}_{\smash{C^0_\beta}}
\cdot\int_{x\in\R^n}\Bigl\vert\md{x-y}^{2-n}-\rho(y)^{2-n}\Bigr\vert
\rho(x)^\beta\d x,
\end{equation*}
and estimating as before shows that $\md{u(y)}\le 
C\nm{f}_{\smash{C^0_\beta}}\rho(y)^{\beta+2}$ when 
$\beta\in(-1-n,-n)$. Thus $u\in C^0_{\beta+2}(\R^n)$ and 
$\nm{u}_{\smash{C^0_{\beta+2}}}\le C\nm{f}_{\smash{C^0_\beta}}$.
The rest of case (b) follows as above.
\end{proof}

Now we extend Theorem \ref{alede1thm} to ALE manifolds.

\begin{thm} Suppose $(X,g)$ is an ALE manifold asymptotic to 
$\R^n/G$ for $n>2$, and\/ $\rho$ a radius function on $X$. 
Let\/ $k\ge 0$ be an integer and\/ $\alpha\in(0,1)$. Then
\begin{itemize}
\item[{\rm(a)}] Let\/ $\beta\in(-n,-2)$. Then there exists $C>0$
such that for each\/ $f\in C^{k,\alpha}_\beta(X)$ there is a unique 
$u\in C^{k+2,\alpha}_{\beta+2}(X)$ with\/ $\Delta u=f$, which
satisfies~$\nm{u}_{\smash{C^{k+2,\alpha}_{\beta+2}}}\le 
C\nm{f}_{\smash{C^{k,\alpha}_\beta}}$.

\item[{\rm(b)}] Let\/ $\beta\in(-1-n,-n)$. Then there exist\/ 
$C_1,C_2>0$ such that for each\/ $f\in C^{k,\alpha}_\beta(X)$ 
there is a unique $u\in C^{k+2,\alpha}_{2-n}(X)$ with\/ 
$\Delta u=f$. Moreover $u=A\rho^{2-n}+v$, where 
\begin{equation}
A={\md{G}\over(n\!-\!2)\,\Omega_{n-1}}\cdot\int_Xf\,\d V_g
\label{adefeq}
\end{equation}
and\/ $v\in C^{k+2,\alpha}_{\beta+2}(X)$ satisfy 
$\md{A}\le C_1\nm{f}_{\smash{C^0_\beta}}$ and\/ 
$\nm{v}_{\smash{C^{k+2,\alpha}_{\beta+2}}}\le 
C_2\nm{f}_{\smash{C^{k,\alpha}_\beta}}$. Here $\Omega_{n-1}$ is 
the volume of the unit sphere ${\mathcal S}^{n-1}$ in~$\R^n$.
\end{itemize}
\label{alede2thm}
\end{thm}

\begin{proof} The theory of weighted H\"older spaces on AE manifolds
and the Laplacian is developed by Chaljub-Simon and Choquet-Bruhat
\cite{ChCh}, who restrict their attention to the case $n=3$. In
particular, they prove part (a) of the Theorem for the case $n=3$,
$k=0$ and $G=\{1\}$, \cite[p.~15-16]{ChCh}. Their proof uses a result 
equivalent to part (a) of Theorem \ref{alede1thm} in the case $n=3$
and $k=0$. By using Theorem \ref{alede1thm} together with the methods
of \cite{ChCh} one can show that Theorem \ref{alede1thm} applies
not only to $\R^n$ with its Euclidean metric, but also to any ALE
manifold $(X,g)$ asymptotic to $\R^n/G$. This proves case (a) of
of the Theorem immediately.

For case (b), let $f\in C^{k,\alpha}_\beta(X)$, and define 
$A$ by \eq{adefeq}. Then by equation \eq{intderhoeq} we have 
$\int_X\bigl[f-\Delta(A\rho^{2-n})\bigr]\d V_g=0$. Also 
$\Delta(\rho^{2-n})\in C^\infty_{-2n}(X)$ by Proposition 
\ref{alede1prop}, and so $f-\Delta(A\rho^{2-n})$ lies in 
$C^{k,\alpha}_\beta(X)$ and has integral zero on $X$. Since 
$\md{f}\le\nm{f}_{\smash{C^0_\beta}}\rho^\beta$ we have 
$\md{A}\le C_1\nm{f}_{\smash{C^0_\beta}}$ for 
$C_1=\int_X\rho^\beta\,\d V_g$, as we have to prove.

Applying case (b) of Theorem \ref{alede1thm} for $X$ to 
$f-\Delta(A\rho^{2-n})$, we see that there is a unique 
$v\in C^{k+2,\alpha}_{\beta+2}(X)$ with 
$\Delta v=f-\Delta(A\rho^{2-n})$, which satisfies
\begin{equation}
\nm{v}_{\smash{C^{k+2,\alpha}_{\beta+2}}}\le 
C\bigl(\nm{f}_{\smash{C^{k,\alpha}_\beta}}+\md{A}\cdot
\nm{\Delta(\rho^{2-n})}_{\smash{C^{k,\alpha}_\beta}}\bigr).
\label{vckaleq}
\end{equation}
Defining $u=A\rho^{2-n}+v$ gives $\Delta u=f$ as we want. 
Clearly $u\in C^{k+2,\alpha}_{2-n}(X)$, and the inequality 
$\nm{v}_{\smash{C^{k+2,\alpha}_{\beta+2}}}\le 
C_2\nm{f}_{\smash{C^{k,\alpha}_\beta}}$ then follows from 
\eq{vckaleq} and the estimate on $\md{A}$ above.
\end{proof}

\section{Exterior forms and de Rham cohomology}

Let $(X,g)$ be an ALE manifold asymptotic to $\R^n/G$. Let 
$H^*(X,\R)$ be the de Rham cohomology of $X$, and $H^*_c(X,\R)$ 
the de Rham cohomology of $X$ {\it with compact support}. That is,
\begin{equation*}
H^k_c(X,\R)={\bigl\{\eta:\text{$\eta$ is a smooth, closed,
compactly-supported $k$-form on $X$}\bigr\}
\over\bigl\{\d\zeta:\text{$\zeta$ is a smooth, compactly-supported 
$(k\!-\!1)$-form on $X$}\bigr\}}.
\end{equation*}
Both $H^k(X,\R)$ and $H^k_c(X,\R)$ are finite-dimensional 
vector spaces. Let us regard $X$ as a compact manifold with 
boundary ${\mathcal S}^{n-1}/G$. Using the long exact sequence
\begin{equation*}
\ldots\rightarrow H^k_c(X,\R)\rightarrow 
H^k(X,\R)\rightarrow H^k({\mathcal S}^{n-1}/G,\R)
\rightarrow H^{k+1}_c(X,\R)\rightarrow\ldots,
\end{equation*}
the de Rham cohomology of ${\mathcal S}^{n-1}/G$, and the fact that
$H^k_c(X,\R)\cong\bigl[H^{n-k}(X,\R)\bigr]^*$ by Poincar\'e duality 
for manifolds with boundary, one can show that
\begin{equation*}
\begin{split}
& H^0(X,\R)=\R,\quad H^0_c(X,\R)=0,\quad
H^n(X,\R)=0,\quad H^n_c(X,\R)=\R,\quad\text{and}\\
& H^k(X,\R)\cong H^k_c(X,\R)\cong
\bigl[H^{n-k}(X,\R)\bigr]^*\cong\bigl[H^{n-k}_c(X,\R)\bigr]^*
\quad\text{for $0<k<n$.}
\end{split}
\end{equation*}

Now the material on weighted H\"older spaces of functions in 
\S 4 generalizes naturally to weighted H\"older spaces of 
$k$-forms on ALE manifolds $(X,g)$, so we may define the 
spaces $C^{l,\alpha}_\beta(\Lambda^kT^*X)$ and 
$C^\infty_\beta(\Lambda^kT^*X)$ in the obvious way. Similarly, 
the results of \S 4 on the Laplacian $\Delta$ on functions 
generalize to results on the Laplacian $\Delta=\d\d^*+\d^*\d$ 
on $k$-forms.

These tools can be used to generalize the ideas of Hodge theory to 
ALE manifolds. In particular, one can prove the following result.

\begin{thm} Let\/ $(X,g)$ be an ALE manifold asymptotic to $\R^n/G$
for $n>2$, and define
\begin{equation*}
\mathcal{H}^k=\bigl\{\eta\in C^\infty_{1-n}(\Lambda^kT^*X):
\d\eta=\d^*\eta=0\bigr\}.
\end{equation*}
Then $\mathcal{H}^0=\mathcal{H}^n=0$, and the map 
$\mathcal{H}^k\rightarrow H^k(X,\R)$ given by 
$\eta\mapsto[\eta]$ induces natural isomorphisms
$\mathcal{H}^k\cong H^k(X,\R)\cong H^k_c(X,\R)$ for $0<k<n$. 
The Hodge star gives an isomorphism $*:\mathcal{H}^k\rightarrow
\mathcal{H}^{n-k}$. Suppose $1-n\le\beta<-n/2$. Then
\begin{equation*}
C^\infty_\beta(\Lambda^kT^*X)=\mathcal{H}^k
\oplus\d\Bigl[C^\infty_{\beta+1}(\Lambda^{k-1}T^*X)\Bigr]
\oplus\d^*\Bigl[C^\infty_{\beta+1}(\Lambda^{k+1}T^*X)\Bigr],
\end{equation*}
where the summands are $L^2$-orthogonal.
\label{alehodgethm}
\end{thm}

This is an analogue of the Hodge Decomposition Theorem and Hodge's
Theorem. For the rest of the section we shall restrict our attention 
to ALE K\"ahler manifolds. If $(X,J,g)$ is an ALE K\"ahler manifold 
then we can define the weighted H\"older spaces of $(p,q)$-forms 
$C^{l,\alpha}_\beta(\Lambda^{p,q}X)$ on $X$ in the obvious way. 
The Laplacian $\Delta$ acts on these spaces by
\begin{equation}
\Delta:C^{l+2,\alpha}_{\beta+2}(\Lambda^{p,q}X)\rightarrow 
C^{l,\alpha}_\beta(\Lambda^{p,q}X).
\end{equation}
They have very similar analytic properties to the weighted H\"older 
spaces of functions on an ALE manifold discussed in~\S 4.

We can use facts about the Laplacian on weighted H\"older spaces 
of $(p,q)$-forms to develop an analogue for ALE K\"ahler manifolds 
of Hodge theory for compact K\"ahler manifolds. 

\begin{thm} Let\/ $(X,J,g)$ be an ALE K\"ahler manifold asymptotic to
$\C^m/G$. Define
\begin{equation}
\mathcal{H}^{p,q}=\bigl\{\eta\in C^\infty_{1-2m}(\Lambda^{p,q}X):
\d\eta=\d^*\eta=0\bigr\}.
\end{equation}
Then $\mathcal{H}^{p,q}$ is finite-dimensional, and the map 
$\mathcal{H}^{p,q}\rightarrow H^{p+q}(X,\C)$ defined by 
$\eta\mapsto[\eta]$ is injective. Define $H^{p,q}(X)$ to be the 
image of this map. Then
\begin{equation}
H^k(X,\C)=\bigoplus_{j=0}^kH^{j,k-j}(X)\qquad\text{for $0<k<2m$.}
\end{equation}
\label{alekhodgethm}
\end{thm}

In fact, if $X$ is a crepant resolution of $\C^m/G$ then 
$H^{p,q}(X)=0$ for~$p\ne q$.

\begin{thm} Let\/ $(X,J,g)$ be an ALE K\"ahler manifold, where $X$ 
is a resolution of\/ $\C^m/G$. Then $H^{2,0}(X)=H^{0,2}(X)=0$, and 
each element of\/ $H^{1,1}(X)$ is represented by a closed, 
compactly-supported\/ $(1,1)$-form on~$X$.
\label{alekh2thm}
\end{thm}

Here is a sketch of the proof of this theorem. Since $X$ is a 
resolution of $\C^m/G$, it can be shown that the homology group
$H_{2m-2}(X,\C)$ is generated by the homology classes of the
exceptional divisors of the resolution. But $H_{2m-2}(X,\C)\cong 
H^2_c(X,\C)$. Thus $H^2_c(X,\C)$ is generated by cohomology classes 
dual to the homology classes $[D]$ of exceptional divisors $D$ in $X$.
If $U$ is any open neighbourhood of $D$ in $X$, then we can find a 
closed $(1,1)$-form supported in $U$ representing the cohomology
class dual to $[D]$. Therefore $H^2_c(X,\C)$ is generated by
cohomology classes represented by closed, compactly-supported
$(1,1)$-forms. It easily follows that $H^{2,0}(X)=H^{0,2}(X)=0$,
and the proof is finished.

Next we prove a version of the Global $\d\d^c$-Lemma for 
ALE K\"ahler manifolds.

\begin{thm} Let\/ $(X,J,g)$ be an ALE K\"ahler manifold asymptotic to 
$\C^m/G$ for some $m>1$, and let\/ $\beta<-m$. Suppose that\/ 
$\eta\in C^\infty_\beta(\Lambda^{1,1}_{\R} X)$ is a closed real 
$(1,1)$-form and\/ $[\eta]=0$ in $H^2(X,\R)$. Then there exists a 
unique real function $u\in C^\infty_{\beta+2}(X)$ with\/~$\eta=\d\d^cu$.
\label{aleddcthm}
\end{thm}

\begin{proof} Let $\omega$ be the K\"ahler form of $g$. Then if 
$u$ is a smooth function on $X$ we have
\begin{equation}
\d\d^cu\wedge\omega^{m-1}=-{\textstyle{1\over m}}\,\Delta u\,\omega^m.
\label{ddcueq}
\end{equation}
Also, if $\zeta$ is a real (1,1)-form on $X$ and 
$\zeta\wedge\omega^{m-1}=0$ it can be shown that
\begin{equation}
\zeta\wedge\omega^{m-2}=-{\textstyle{1\over 2}}(m-2)!\,*\zeta
\quad\text{and}\quad
\zeta\wedge\zeta\wedge\omega^{m-2}=
-{\textstyle{1\over 2}}(m-2)!\,\ms{\zeta}\d V_g,
\label{zetomeq}
\end{equation}
where $*$ is the Hodge star and $\d V_g$ the volume form of $g$.
Equations \eq{ddcueq} and \eq{zetomeq} hold on any K\"ahler 
manifold of dimension~$m$.
 
Define a function $f$ on $X$ by 
$\eta\wedge\omega^{m-1}=-{1\over m}f\,\omega^m$. Since 
$\eta\in C^\infty_\beta(\Lambda^{1,1}_{\R} X)$, it follows that 
$f\in C^\infty_\beta(X)$. Now suppose for simplicity that 
$-2m<\beta<-m$. Then by part (a) of Theorem \ref{alede2thm} 
there exists a unique function $u\in C^\infty_{\beta+2}(X)$ 
with $\Delta u=f$. Set $\zeta=\eta-\d\d^cu$, which is an exact 
2-form in $C^\infty_\beta(\Lambda^{1,1}_{\R} X)$. As $\beta<-m$ 
we can use the last part of Theorem \ref{alehodgethm} to prove that 
$\zeta=\d\theta$, for some~$\theta\in C^\infty_{\beta+1}(T^*X)$. 

By \eq{ddcueq} we have $\zeta\wedge\omega^{m-1}
=-{1\over m}(f-\Delta u)\,\omega^m=0$, so \eq{zetomeq} gives
\begin{equation}
\d\bigl[\theta\wedge\zeta\wedge\omega^{m-2}\bigr]=
\zeta\wedge\zeta\wedge\omega^{m-2}=
-{\textstyle{1\over 2}}(m\!-\!2)!\,\ms{\zeta}\d V_g.
\label{dgazetaeq}
\end{equation}
Let $\rho$ be a radius function on $X$, and define 
$S_R=\bigl\{x\in X:\rho(x)\le R\bigr\}$ for $R>1$. Integrating
\eq{dgazetaeq} over $S_R$ and using Stokes' Theorem gives that
\begin{equation}
-{\textstyle{1\over 2}}(m-2)!\cdot\int_{S_R}\ms{\zeta}\d V_g=
\int_{\partial S_R}\theta\wedge\zeta\wedge\omega^{m-2}.
\label{intsrzeta2eq}
\end{equation}
But for large $R$ we have $\theta=O(R^{\beta+1})$, 
$\zeta=O(R^\beta)$ and $\omega=O(1)$ on $\partial S_R$, and 
$\vol(\partial S_R)=O(R^{2m-1})$. Thus the r.h.s.~of 
\eq{intsrzeta2eq} is $O(R^{2\beta+2m})$. As $\beta<-m$, 
taking the limit as $R\rightarrow\infty$ shows that 
$\int_X\ms{\zeta}\d V_g=0$, and so $\zeta=0$ on $X$. 
Thus $\eta=\d\d^cu$, as we have to prove.

We have proved the theorem assuming that $-2m<\beta<-m$, 
but we wish to prove it for all $\beta<-m$. If $\beta\le 2m$ 
and $\eta\in C^\infty_\beta(\Lambda^{1,1}_{\R} X)$ then 
$\eta\in C^\infty_\gamma(\Lambda^{1,1}_{\R} X)$ for any 
$\gamma$ with $-2m<\gamma<-m$, and so from above we have 
$\eta=\d\d^cu$ for some unique $u$ in $C^\infty_{\gamma+2}(X)$. 
However, if $u\in C^\infty_{\gamma+2}(X)$ and
$\d\d^cu\in C^\infty_\beta(\Lambda^{1,1}_{\R} X)$, one can 
show that $u\in C^\infty_{\beta+2}(X)$ as we want. This is 
because $\d\d^cu$ is a stronger derivative of $u$ than 
$\Delta u$ is, and contains more information.
\end{proof}

Finally, we show we can modify any ALE K\"ahler metric to 
be flat outside a compact set.

\begin{prop} Let\/ $\C^m/G$ have an isolated singularity at\/ 
$0$ for some $m>1$, let\/ $(X,\pi)$ be a resolution of\/ 
$\C^m/G$ that admits ALE K\"ahler metrics, and let\/ $\rho$ 
be a radius function on $X$. Then in each K\"ahler class 
there exists an ALE K\"ahler metric $\hat g$ on $X$ such that\/ 
$\hat g=\pi^*(h)$ on the subset\/ $\bigl\{x\in X:\rho(x)>R\bigr\}$, 
where $h$ is the Hermitian metric on $\C^m/G$ and\/ $R>0$ is 
a constant.
\label{alekprop}
\end{prop}

\begin{proof} Let $g$ be an ALE K\"ahler metric on $X$, with 
K\"ahler form $\omega$. By Theorems \ref{alekhodgethm} and 
\ref{alekh2thm} there exists a closed, compactly-supported, 
real (1,1)-form $\theta$ on $X$ with $[\theta]=[\omega]$ in 
$H^2(X,\R)$. Define $\eta=\omega-\d\d^c(\rho^2)-\theta$. Then 
$\eta$ is an exact real (1,1)-form on $X$. Now the K\"ahler 
form of $h$ on $\C^m/G$ is $\omega_0=\d\d^c(r^2)$. So from 
the definition of ALE K\"ahler metric we see that 
$\omega-\d\d^c(\rho^2)\in C^\infty_{-2m}(\Lambda^{1,1}_{\R} X)$, 
and therefore $\eta\in C^\infty_{-2m}(\Lambda^{1,1}_{\R} X)$ 
as $\theta$ has compact support. Thus by Theorem \ref{aleddcthm} 
there is a unique real function $u\in C^\infty_{2-2m}(X)$ with 
$\eta=\d\d^cu$, and we have~$\omega=\theta+\d\d^c(\rho^2)+\d\d^cu$.

Let $\mu:\R\rightarrow[0,1]$ be a smooth function with $\mu(t)=1$ 
for $t\le -1$ and $\mu(t)=0$ for $t\ge 0$. For each $R>0$ define
a closed (1,1)-form $\omega_R$ by
\begin{equation}
\omega_R=\theta+\d\d^c(\rho^2)+\d\d^c\bigl[\mu(\rho-R)\cdot u\bigr].
\end{equation}
Then $\omega_R=\omega$ wherever $\rho<R-1$, and $\omega_R=\d\d^c(\rho^2)$
wherever $\rho>R$ and outside the support of $\theta$. It is 
easy to show that $\omega_R$ is a positive (1,1)-form for large 
$R$, which therefore defines a K\"ahler metric $g_R$ on $X$. 
Define $\hat g$ to be $g_R$ for some $R$ sufficiently large that
$\omega_R$ is positive, $\rho\le R$ on the support of $\theta$
and $R\ge 2$. Then $\hat g$ is an ALE K\"ahler metric in the K\"ahler 
class of $g$, and where $\rho>R$ we have $\hat g=\pi^*(h)$, since the
K\"ahler form of $\hat g$ is $\d\d^c(\rho^2)$, the K\"ahler form of 
$h$ is $\d\d^c(r^2)$, and $\rho=\pi^*(r)$ as~$\rho>R\ge 2$.
\end{proof}

\section{The Calabi conjecture for ALE manifolds}

We can now state the following version of the Calabi 
conjecture for ALE K\"ahler manifolds.
\medskip

\noindent{\bf The Calabi conjecture for ALE manifolds} 
{\it Suppose that\/ $(X,J,g)$ is an ALE K\"ahler manifold of 
dimension $m$ asymptotic to $\C^m/G$ for some $m>1$, with 
K\"ahler form $\omega$, and that\/ $\rho$ is a radius function 
on $X$. Then
\begin{itemize}
\item[{\rm(a)}] Let\/ $\beta\in(-2m,-2)$. Then for 
each\/ $f\in C^\infty_\beta(X)$ there is a unique 
$\phi\in C^\infty_{\beta+2}(X)$ such that\/ 
$\omega+\d\d^c\phi$ is a positive $(1,1)$-form and\/ 
$(\omega+\d\d^c\phi)^m={\rm e}^f\omega^m$ on~$X$.
\item[{\rm(b)}] Let\/ $\beta\in(-1-2m,-2m)$. 
Then for each\/ $f\in C^\infty_\beta(X)$ there is a 
unique $\phi\in C^\infty_{2-2m}(X)$ such that\/ 
$\omega+\d\d^c\phi$ is a positive $(1,1)$-form and\/ 
$(\omega+\d\d^c\phi)^m={\rm e}^f\omega^m$ on $X$. 
Moreover we can write $\phi=A\rho^{2-2m}+\psi$, 
where $\psi\in C^\infty_{\beta+2}(X)$ and
\begin{equation}
A={\md{G}\over(m-1)\Omega_{2m-1}}\cdot\int_X(1-{\rm e}^f)\d V_g.
\label{aleccadefeq}
\end{equation}
Here $\Omega_{2m-1}$ is the volume of the unit sphere 
${\mathcal S}^{2m-1}$ in~$\C^m$.
\end{itemize}}
\medskip

It is easy to rewrite this in terms of the existence of ALE 
K\"ahler metrics with prescribed Ricci curvature, as in the 
original Calabi conjecture. The two cases (a) $\beta\in(-2m,-2)$ 
and (b) $\beta\in(-1-2m,-2m)$ come from Theorem \ref{alede2thm}. 
By combining the method of Yau's proof \cite{Yau} of the 
Calabi conjecture with the ideas of \S 4 on analysis
on ALE manifolds, we can prove the Calabi conjecture for
ALE manifolds.

The conjecture will be proved in \cite[\S 8.5--\S 8.6]{Joyc3}, 
and we give only a sketch of the proof of part (a) here. We use 
the {\it continuity method}. Suppose $\beta\in(-2m,-2)$. Fix 
$f\in C^{3,\alpha}_\beta(X)$, and define $S$ to be the set of all 
$t\in[0,1]$ for which there exists $\phi\in C^{5,\alpha}_{\beta+2}(X)$ 
such that $\omega+\d\d^c\phi$ is a positive (1,1)-form and 
$(\omega+\d\d^c\phi)^m={\rm e}^{tf}\omega^m$ on~$X$. 

Clearly $0\in S$, taking $\phi=0$. We prove that $S$ is 
both {\it open} and {\it closed} in $[0,1]$. Thus $S=[0,1]$ 
as $[0,1]$ is connected, so $1\in S$, and there exists 
$\phi\in C^{5,\alpha}_{\beta+2}(X)$ with $\omega+\d\d^c\phi$ 
positive and $(\omega+\d\d^c\phi)^m={\rm e}^f\omega^m$ on $X$.
We then use Theorem \ref{alede2thm} to show that if 
$f\in C^\infty_\beta(X)$ then $\phi\in C^\infty_{\beta+2}(X)$, 
and this completes the proof. 

To prove that $S$ is open, we fix $t\in S$ and show that $S$
contains a small neighbourhood of $t$ by considering the {\it 
linearization} of the equation at $t$. This linearization
turns out to involve the Laplacian of the metric with K\"ahler 
form $\omega+\d\d^c\phi$, and part (a) of Theorem \ref{alede2thm} 
gives us what we need.

To prove that $S$ is closed, we take a sequence 
$\{t_j\}_{j=0}^\infty$ in $S$ such that $t_j\rightarrow t\in [0,1]$ as 
$j\rightarrow\infty$. Let $\{\phi_j\}_{j=0}^\infty$ be the sequence of 
solutions to $(\omega+\d\d^c\phi_j)^m={\rm e}^{t_jf}\omega^m$. 
Then $\phi_j$ converges to some $\phi\in C^{5,\alpha}_{\beta+2}(X)$ as 
$j\rightarrow\infty$ with $(\omega+\d\d^c\phi)^m={\rm e}^{tf}\omega^m$, 
and thus $t\in S$. Therefore $S$ contains its limit points, and is closed.

The difficult part in showing $S$ closed is finding an {\it a 
priori estimate} for $\phi_j$ in $C^{5,\alpha}_{\beta+2}(X)$. 
To do this we first follow Yau's proof to get an a priori estimate 
in $C^{5,\alpha}(X)$. Then we use a `weighted' version of Yau's 
method to estimate $\phi_j$ in $C^0_\delta(X)$ for some small 
$\delta<0$. This can be improved to $C^{5,\alpha}_\delta(X)$, and 
then to $C^{5,\alpha}_{\beta+2}(X)$ by a kind of induction, 
decreasing $\delta$ step by step until~$\delta=\beta+2$.

This concludes our treatment of the Calabi conjecture 
for ALE manifolds, and we are now ready to prove 
Theorems \ref{alerfthm} and~\ref{alesumthm}.

\subsection{The proof of Theorem \ref{alerfthm}}

Let $X$ be a crepant resolution of $\C^m/G$, where $G$ acts freely
on $\C^m\setminus\{0\}$. By Proposition \ref{alekprop}, in each 
K\"ahler class of ALE K\"ahler metrics on $X$ we can choose a 
metric $\hat g$ with $\hat g=\pi^*(h)$ wherever $\rho>R\ge 2$, 
where $h$ is the Euclidean metric on $\C^m/G$. Let $\hat\omega$ 
be the K\"ahler form and $\eta$ the Ricci form of $\hat g$. Then 
$\eta$ is closed and $[\eta]=2\pi\,c_1(X)$ in $H^2(X,\R)$. But 
$c_1(X)=0$ as $X$ is a crepant resolution, so $[\eta]=0$ in 
$H^2(X,\R)$. Also, $\eta=0$ wherever $\rho>R$, since there 
$\hat g=\pi^*(h)$ and $h$ is flat. 

Thus $\eta$ is a closed, compactly-supported $(1,1)$-form on $X$ 
with $[\eta]=0$ in $H^2(X,\R)$, and by Theorem \ref{aleddcthm} 
there exists a unique function $f\in C^\infty_\beta(X)$ for each 
$\beta<0$ with $\eta={1\over 2}\d\d^cf$. In fact $f=0$ wherever 
$\rho>R$, so $f$ is compactly supported. The Calabi conjecture 
for ALE manifolds holds by \cite[\S 8.5--\S 8.6]{Joyc3}. Part 
(b) of the conjecture shows that there exists a unique function 
$\phi=A\rho^{2-2m}+\psi$ where $A$ is given by \eq{aleccadefeq} 
and $\psi\in C^\infty_{\beta+2}(X)$ for $\beta\in(-1-2m,-2m)$, 
such that $\omega=\hat\omega+\d\d^c\phi$ is a positive 
(1,1)-form and~$\omega^m={\rm e}^f\hat\omega^m$. 

Let $g$ be the K\"ahler metric on $X$ with K\"ahler form $\omega$. 
Then since the Ricci form of $\hat g$ is ${1\over 2}\d\d^cf$ it 
follows by standard properties of the Ricci form that $g$ has Ricci 
form zero, and is Ricci-flat. On $\bigl\{z\in\C^m/G:r(z)>R\bigr\}$ 
we have $\pi_*(\hat g)=h$, so that $\pi_*(\hat\omega)=\omega_0$, 
and $\pi_*(\rho)=r$. Thus defining $\chi=\pi_*(\psi)$ gives 
\eq{pisomeq}. Since $\psi\in C^\infty_{\beta+2}(X)$ for 
$\beta\in(-1\!-2m,-2m)$, putting $\gamma=\beta+2$ we see that
$\nabla^k\chi=O(r^{\gamma-k})$ for $k=0,1,2,\ldots$ and 
$\gamma\in(1-2m,2-2m)$, as we have to prove. 

{}From \eq{pisomeq} we see that $\nabla^k(\pi_*(g)-h)=O(r^{-2m-k})$ 
for $k\ge 0$, and thus $g$ is an ALE K\"ahler metric by Definition 
\ref{alekdef}. Also, $g$ is unique in its K\"ahler class of ALE 
metrics because $\phi$ is unique. It only remains to prove that 
$A<0$. We can do this by giving an explicit expression for $A$. 
Let $\zeta$ be the unique element of $\mathcal{H}^{1,1}$ with 
$[\zeta]=[\omega]$. Then a calculation shows that
\begin{equation}
A=-\,{\md{G}\over 2m(m\!-\!1)^2\Omega_{2m-1}}\int_X\ms{\zeta}\d V_g,
\label{aconsteq}
\end{equation}
where $\Omega_{2m-1}$ is the volume of the unit sphere 
${\mathcal S}^{2m-1}$ in $\C^m$. Now $[\omega]\ne 0$ as this 
is outside the K\"ahler cone, so $\zeta\ne 0$ and $A$ is 
negative. This completes the proof of Theorem~\ref{alerfthm}.

\subsection{The proof of Theorem \ref{alesumthm}}

First we show that $\C^m/G$ has no crepant resolutions when
$m>2$ and~$G\subset{\rm Sp}(m/2)$.

\begin{prop} Suppose that\/ $m>2$ is even and that\/ $G$ is a 
nontrivial finite subgroup of\/ ${\rm Sp}(m/2)$ which acts 
freely on $\C^m\setminus\{0\}$. Then $\C^m/G$ is a terminal 
singularity, and admits no crepant resolutions.
\label{spm2prop}
\end{prop}

\begin{proof} Let $\gamma\ne 1$ in $G$. Then there are coordinates 
$(z^1,\ldots,z^m)$ on $\C^m$ in which $\gamma$ acts by
\begin{equation}
\bigl(z^1,\ldots,z^m\bigr)\,{\buildrel\gamma\over\longmapsto}\,
\bigl({\rm e}^{2\pi ia_1}z^1,\ldots,{\rm e}^{2\pi ia_m}z^m\bigr).
\label{gajacteq}
\end{equation}
As $\gamma$ acts freely on $\C^m\setminus\{0\}$ we can take
$a_j\in (0,1)$ for $j=1,\ldots,m$. Since $G\subset{\rm Sp}(m/2)$ 
we know that $\gamma$ preserves a complex symplectic form on 
$\C^m$, and we can choose $(z^1,\ldots,z^m)$ so that this form 
is $\d z^1\wedge\d z^2+\cdots+\d z^{m-1}\wedge\d z^m$. Thus 
\eq{gajacteq} gives ${\rm e}^{2\pi ia_{2j-1}}{\rm e}^{2\pi ia_{2j}}=1$ 
for $j=1,\ldots,m/2$. But as $a_{2j-1},a_{2j}\in (0,1)$ this implies 
that $a_{2j-1}+a_{2j}=1$ for~$j=1,\ldots,m/2$.

Therefore $a_1+\cdots+a_m=m/2>1$ for all $\gamma\ne 1$ in $G$.
So by Reid \cite[\S 4]{Reid} it follows that $\C^m/G$ is a 
{\it terminal singularity}, as defined in~\cite[p.~347]{Reid}.
Terminal singularities are essentially singularities which have
no crepant partial resolutions. To be more precise, a crepant
resolution of a terminal singularity has no exceptional divisors.
Thus, if $X$ is a crepant resolution of $\C^m/G$ then $b_{2m-2}(X)=0$.
By Poincar\'e duality for manifolds with boundary we see that
$b_2(X)=0$, which is a contradiction, as $X$ must contain a
complex curve. So $\C^m/G$ has no crepant resolutions.
\end{proof}

We now prove Theorem \ref{alesumthm}. Let $X$ be a crepant
resolution of $\C^m/G$, where $G$ is nontrivial and acts
freely on $\C^m\setminus\{0\}$, and let $g$ be a Ricci-flat ALE
K\"ahler metric on $X$. As $X$ is simply-connected, by general
facts about holonomy groups we know that $\Hol(g)$ is a 
connected Lie subgroup of ${\rm SU}(m)$. Since $g$ is Ricci-flat 
it is nonsymmetric. Also $(X,g)$ is not a Riemannian product, 
because it is asymptotic to $\C^m/G$, which is not a product. 
Thus $g$ is irreducible. 

Therefore we may apply Berger's classification of Riemannian 
holonomy groups \cite[\S 10]{Sala}. The only two possibilities 
are $\Hol(g)={\rm SU}(m)$ or $\Hol(g)={\rm Sp}(m/2)$. When 
$m=2$ the two groups coincide, so suppose $m>2$. The holonomy 
of the Euclidean metric $h$ on $\C^m/G$ is $G\subset{\rm SU}(m)$. 
Since $g$ is asymptotic to $h$ one can show that 
$G\subset\Hol(g)\subseteq{\rm SU}(m)$. Hence, if 
$\Hol(g)={\rm Sp}(m/2)$ then $G\subset{\rm Sp}(m/2)$. 
But Proposition \ref{spm2prop} then shows that $\C^m/G$ 
admits no crepant resolutions, a contradiction. So 
$\Hol(g)\ne{\rm Sp}(m/2)$, and thus $\Hol(g)={\rm SU}(m)$, 
which completes the proof.

The essential point in this proof is that there do not exist
ALE manifolds with holonomy ${\rm Sp}(m/2)$ for $m>2$. One 
can also show this using Schlessinger's Rigidity Theorem 
\cite{Schl}, and properties of hyperk\"ahler manifolds.

\end{document}